\documentclass[11pt]{amsart}
\usepackage{amsmath,amssymb,amscd,amsthm,amsxtra}

\allowdisplaybreaks

\newtheorem{theorem}{Theorem}

\newtheorem{lemma}[theorem]{Lemma}

\def\R{\mathbb{R}}
\def\T{\mathbb{T}}

\newcommand{\norm}[2]{\left\| #1 \right\|_{#2}}
\newcommand{\mdl}[1]{\left| #1 \right|}

\begin{document}

\title{Perfect dyadic operators: weighted $T(1)$ theorem and two weight estimates.}

\author[O. Beznosova]{Oleksandra Beznosova}

\address{Oleksandra Beznosova\\
Department of Mathematics\\
University of Alabama\\ 345 Gordon Palmer Hall\\ Tuscaloosa, AL 35487}
\email{ovbeznosova@ua.edu}

\thanks{The author was supported by the University of Alabama RGC grant.}


\begin{abstract}
Perfect dyadic operators were first introduced in \cite{AHMTT}, where a local $T(b)$ theorem was proved for such operators. In \cite{AY} it was shown that for every singular integral operator $T$ with locally bounded kernel on $\R^n \times \R^n$ there exists a perfect dyadic operator $\T$ such that $T -\T$ is bounded on $L^p (dx)$ for all $1<p<\infty$.

In this paper we show a decomposition of perfect dyadic operators on real line into four well known operators: two selfadjoint operators, paraproduct and its adjoint. Based on this decomposition we prove a sharp weighted version  of the $T(1)$ theorem for such operators, which implies $A_2$ conjecture for such operators with constant which only depends on $\|T(1)\|_{BMO^d}$, $\|T^*(1)\|_{BMO^d}$ and the constant in testing conditions for $T$. Moreover, the constant depends on these parameters at most linearly. In this paper we also obtain sufficient conditions for the two weight boundedness for a perfect dyadic operator and simplify these conditions under additional assumptions that weights are in the Muckenhoupt class $A_\infty^d$.
\end{abstract}



\maketitle

\section{Introduction}


A perfect dyadic operator $T$ is defined by:
\begin{equation} \nonumber
\label{formula 1 T} Tf \; := \; \int K(x,y) \, f(y) \, dy \;\; \text{for} \;\; x \notin supp\, f,
\end{equation}
where kernel $K$ satisfies the following conditions:\\
standard size condition:
\begin{equation}
\label{formula 2 cond}
\mdl{K(x,y)} \; \leq \; \frac{1}{\mdl{x-y}}
\end{equation}
and perfect cancellation condition:
\begin{equation}
\label{formula 3 cond}
\mdl{K(x,y)-K(x,y^\prime)} +
\mdl{K(x,y)-K(x^\prime,y)} \; = \; 0
\end{equation}
whenever $x,x^\prime \in I \in D$, $y,y^\prime \in J \in D$, $I \cap
J = \varnothing$. Where $D$ is a set of dyadic intervals on the real line $D:= \{I = [2^j k, 2^j (k+1)): \; k, j \in \mathbb{Z}\}$.

Perfect dyadic operators appeared in the context of the local $T(b)$ theorems (see \cite{AHMTT}) and later were extended to spaces of homogeneous type (see \cite{AR}, \cite{LV}). Perfect dyadic operators are essentially Calder\'{o}n--Zygmund singular integrals with singularity adapted to the fixed dyadic grid. Their main property is that any function supported on a dyadic cube with zero mean is mapped into the function supported on the same cube. In \cite{AY} it was shown that in $L^2(dx)$ perfect dyadic operators are a good approximation of Calder\'{o}n--Zygmund singular integral operators. 

In this paper we stick to the real line as a model case. In Section 2 we derive a very useful decomposition of the perfect dyadic operator into four well known operators: two selfadjoint operators, dyadic paraproduct and its adjoint. In particular, we show that a perfect dyadic operator $T$ and its adjoint formally can be written as

\begin{eqnarray}
\label{e:decompT}\nonumber Tf(x) &=& \frac{1}{4}\sum_{I\in
D}\left( K_I^+ + K_I^- \right) \mdl{I} \, m_If \, \chi_I -
 \frac{1}{4}\sum_{I\in
D}\left( K_I^+ - K_I^- \right)\sqrt{\mdl{I}}
\left\langle f;h_I \right\rangle \chi_I\\
&& +  \frac{1}{4}\sum_{I\in
D}\left( K_I^+ - K_I^- \right)
\mdl{I}^{3/2} \, m_If \, h_I -  \frac{1}{4}\sum_{I\in
D}\left( K_I^+
+ K_I^- \right)\mdl{I}\left\langle f;h_I \right\rangle h_I,
\end{eqnarray}
\begin{eqnarray}
\label{e:decompT*}\nonumber T^*f(x) &=& \frac{1}{4}\sum_{I\in D}\left( K_I^+ + K_I^- \right)
\mdl{I} \,m_If\, \chi_I(x) +  \frac{1}{4}\sum_{I\in
D}\left( K_I^+ -
K_I^- \right)\sqrt{\mdl{I}} \left\langle f;h_I \right\rangle
\chi_I(x)\\
&& -  \frac{1}{4}\sum_{I\in
D}\left( K_I^+ - K_I^- \right)
\mdl{I}^{3/2} m_If \, h_I(x) -  \frac{1}{4}\sum_{I\in
D}\left( K_I^+
+ K_I^- \right) \mdl{I} \left\langle f;h_I \right\rangle h_I(x).
\end{eqnarray}

For a given dyadic interval $I$, $I^+$ and $I^-$ are its left and right halves. Coefficients $K_I^+$ and $K_I^-$ are defined to be values of the kernel $K$ of the perfect dyadic operator $T$ on the dyadic cubes $I^+ \times I^-$ and $I^- \times I^+$ respectively. The notation $m_I f$ stands for the average of the function $f$ over the interval $I$, $m_I f := \frac{1}{|I|} \int_I f$.

In Section 3 using our decomposition we prove the $T(1)$ theorem and show several useful estimates on the kernel. We prove the following theorem.

\begin{theorem} \label{t:T(1)}
Let $T$ be a perfect dyadic operator that satisfies:

(i) BMO conditions: $\|T(1)\|_{BMO^d}$ and $\|T^*(1)\|_{BMO^d}$ are bounded by a finite numerical constant $Q$;

(ii) Testing conditions: $\langle T h_I; h_I \rangle$ are uniformly bounded by $Q$ for every dyadic interval $I\in D$.

Then $T$ is bounded on $L^2$.

Moreover, $T$ accepts decomposition (\ref{e:decompT}) with coefficients $K_I^\pm$ satisfying

\noindent size conditions:
\begin{equation}
\label{formula 19 size cond} \forall J\in D \;\;\;\;\; \mdl{K_J^+ +
K_J^-} \mdl{I} \; \leq \; 2 Q,
\end{equation}
$T(1)$ conditions:
\begin{equation}
\label{formula 20 T(1) cond} \forall J\in D \;\;\;\;\;
\frac{1}{\mdl{J}}\sum_{I\in D(J)}\left( K_I^+ - K_I^- \right)^2
\mdl{I}^3 \; \leq \; 16 Q
\end{equation}
and testing conditions:
\begin{equation}
\label{formula 21 test cond} \forall J\in D \;\;\;\;\;
\mdl{\frac{1}{\mdl{J}}\sum_{I\in D(J)}\left( K_I^+ + K_I^-
\right)\mdl{I}^2} \; \leq \; 8 Q.
\end{equation}
\end{theorem}

In Section 4 we `lift' the above $T(1)$ theorem to the one weight case, $L^2(w)$, and provide an elementary proof of the $A_2$ conjecture for such operators; in particular we trace how the weighted norm of perfect dyadic operator depends on the dyadic $BMO^d$ and testing constants of the perfect dyadic operator. We prove the following theorem.

\begin{theorem} 
Let $T$ be a perfect dyadic operator on $\R$ such that 
$$\|T(1)\|_{BMO^d} \leqslant Q \;\;\; \text{and} \;\;\; \|T^*(1)\|_{BMO^d} \leqslant Q $$
and
$$ \forall I\in D \;\;\; \langle T h_I , h_I \rangle \leqslant Q $$.
Then $T$ is bounded on $L^2(w)$ and 
\begin{equation}\label{e:A2}
\|T f\|_{L^2(w)} \leqslant C Q [w]_{A_2^d}\|f\|_{L^2(w)},
\end{equation}
with some constant $C$ independent of the operator $T$.
\end{theorem}

In particular, by \cite{DGPP}, this implies that under assumptions of the theorem a perfect dyadic operator $T$ is also bounded in $L^p (w)$ for all $w\in A_p^d$, its norms are bounded by $C [w]_{A_p^d}^{\max\{1, \frac{1}{p}\}}$, and dependence $C(Q)$ can be traced as well. 


To the best of our knowledge, such sharp weighted version of $T(1)$ theorem is new. In particular it is interesting that the constant in (\ref{e:A2}) depends only on the $BMO^d$ norms of $T(1)$ and $T^*(1)$ and on the constant in testing conditions. Moreover, the dependence is at most linear. 

In Section 5 we go even further and give sufficient conditions for the two weight boundedness for a single perfect dyadic operator. Currently there are two schools of thought regarding the two weight problem. First, given one operator $T$ find necessary and sufficient conditions
on the weights to ensure boundedness of the operator on the appropriate spaces. Second, given a family of operators
find necessary and sufficient conditions on the weights to ensure boundedness of the family of operators. In both cases we are mostly interested in Calder\'on--Zygmund singular integral operators.
In the first case, the conditions are usually ``testing conditions" obtained from checking boundedness of the given operator on
a collection of test functions. In the second case, the conditions are more ``geometric", meaning they seem to only involve the weights and
not the operators, such as Carleson conditions or bilinear embedding conditions, Muckenhoupt $A_2$ type conditions or bumped conditions.
Operators of interest are the maximal function \cite{S1, Moe, PzR, V}, fractional and Poisson integrals \cite{S2,Cr}, the Hilbert transform \cite{CS1, CS2, KP, NazarovTreilVolberg:99, LSSU, L3}
and general Calder\'on--Zygmund singular integral operators and their commutators \cite{CrRV, CrMoe,CrMPz2, HoLWic,NRTV}, the square functions \cite{ LLi,LLi2}, paraproducts and their dyadic counterparts. Necessary and sufficient conditions are only known for the maximal function, fractional and Poisson integrals \cite{S1}, square functions \cite{LLi} and the Hilbert transform \cite{L3,LSSU}, and among the dyadic operators for the martingale transform, the dyadic square functions, positive and well localized dyadic operators  \cite{Wil87, NazarovTreilVolberg:99,NTV3, T, Ha, HaHLi, HL, LSU2, Ta, Vu}.


We prove the following sufficient conditions for the two weight boundedness of the perfect dyadic operators.

\begin{theorem} \label{t:2w-1}
Let $(v,u)$ be a pair of measurable functions, such that $u$ and $v^{-1}$, the reciprocal of $v$, are weights on $\R$ such that
\begin{equation}\label{e:uv-1}
(v,u)\in A_2^d \,,\text{i.e.} \;\; [v,u]_{A_2^d} := \sup_{I\in D} m_I (v^{-1}) m_I u <\infty,
\end{equation}
\begin{equation}\label{e:uv-2}
\forall J \in D \;\;\; \frac{1}{|J|} \sum_{I \in D(J)} (\Delta_I u)^2 m_I (v^{-1}) |I| \leqslant C m_I u,
\end{equation}
\begin{equation}\label{e:uv-3}
\forall J \in D \;\;\; \frac{1}{|J|} \sum_{I \in D(J)} (\Delta_I (v^{-1}))^2 m_I u |I| \leqslant C m_I (v^{-1}),
\end{equation}
and operator $T_0$ is bounded from $L^2(v)$ in $L^2(u)$, i.e.
\begin{equation}\label{e:uv-4}
\forall J \in D \;\;\; \frac{1}{|J|} \int_J \left(\sum_{I \in D(J)}  \left|\frac{\Delta_I u \Delta_I (v^{-1})}{m_Iu } \right| \chi_I(x) \right)^2 u(x) dx \leqslant C m_J (v^{-1})
\end{equation}
and
\begin{equation}\label{e:uv-5}
\forall J \in D \;\;\; \frac{1}{|J|} \int_J \left(\sum_{I \in D(J)}   \left|\frac{\Delta_I u \Delta_I (v^{-1})}{ m_I(v^{-1})}\right| \chi_I(x) \right)^2 v^{-1}(x) dx \leqslant C m_J u.
\end{equation}
Let $T$ be a perfect dyadic operator with $\|T(1)\|_{BMO^d}, \|T^*(1)\|_{BMO^d} \leqslant Q$ and for every $I \in D$ we have $\langle T h_I; h_I \rangle \leqslant Q$.

Then $T$ is bounded from $L^2 (v)$ to $L^2(u)$ whenever for any dyadic interval $J \in D$

\begin{equation}\label{e:Tuv-1}
\frac{1}{|J|} \sum_{I\in D(J)} (K_I^+ + K_I^-) |I|^2  m_I u \leqslant C m_J u,
\end{equation}
\begin{equation}\label{e:Tuv-2}
\frac{1}{|J|} \sum_{I\in D(J)} (K_I^+ + K_I^-) |I|^2 m_I(v^{-1}) \leqslant C m_J (v^{-1}),
\end{equation}
\begin{equation}\label{e:Tuv-3}
\frac{1}{|J|}\sum_{I\in {D}(J)} \frac{(K_I^+ - K_I^-)^2 |I|^{3}}{m_Iu} \leq C m_J (v^{-1})
\end{equation}
and
\begin{equation}\label{e:Tuv-4}
\frac{1}{|J|}\sum_{I\in {D}(J)} \frac{(K_I^+ - K_I^-)^2 |I|^{3}}{m_I (v^{-1})} \leq C m_J u. 
\end{equation}
\end{theorem}


When both weights are in the Muckenhoupt class $A_\infty^d$ we simplify conditions of Theorem \ref{t:2w-1} to the following.

\begin{theorem}
Let $(v,u)$ be a pair of measurable functions, such that $u$ and $v^{-1}$, the reciprocal of $v$, are $A_\infty^d$ weights on $\R$ such that
\begin{equation}\label{e:uv-1}
(v,u)\in A_2^d \,,\text{i.e.} \;\; [v,u]_{A_2^d} := \sup_{I\in D} m_I (v^{-1}) m_I u <\infty,
\end{equation}
and operator $T_0$ is bounded from $L^2(v)$ in $L^2(u)$, i.e.
\begin{equation}\label{e:uv-4}
\forall J \in D \;\;\; \frac{1}{|J|} \int_J \left(\sum_{I \in D(J)}  \left|\frac{\Delta_I u \Delta_I (v^{-1})}{m_Iu } \right| \chi_I(x) \right)^2 u(x) dx \leqslant C m_J (v^{-1})
\end{equation}
and
\begin{equation}\label{e:uv-5}
\forall J \in D \;\;\; \frac{1}{|J|} \int_J \left(\sum_{I \in D(J)}   \left|\frac{\Delta_I u \Delta_I (v^{-1})}{ m_I(v^{-1})}\right| \chi_I(x) \right)^2 v^{-1}(x) dx \leqslant C m_J u.
\end{equation}
Let $T$ be a perfect dyadic operator with $\|T(1)\|_{BMO^d}, \|T^*(1)\|_{BMO^d} \leqslant Q$ and for every $I \in D$ we have $\langle T h_I; h_I \rangle \leqslant Q$.

Then $T$ is bounded from $L^2 (v)$ to $L^2(u)$ whenever for any dyadic interval $J \in D$

\begin{equation}\label{e:Tuv-3}
\frac{1}{|J|}\sum_{I\in {D}(J)} \frac{(K_I^+ - K_I^-)^2 |I|^{3}}{m_Iu} \leq C m_J (v^{-1})
\end{equation}
and
\begin{equation}\label{e:Tuv-4}
\frac{1}{|J|}\sum_{I\in {D}(J)} \frac{(K_I^+ - K_I^-)^2 |I|^{3}}{m_I (v^{-1})} \leq C m_J u. 
\end{equation}
\end{theorem}

The main novelty of this paper is the decomposition, very strong quantitative estimates in the $A_2$ conjecture and sufficient conditions for the two weight boundedness of perfect dyadic operators. 


\section{Perfect dyadic operator on $\R$: definition and formal decomposition.}

Let $T$ be a perfect dyadic singular integral operator, i.e an operator defined by:
\begin{equation}
\label{formula 1 T} Tf \; := \; \int K(x,y) \, f(y) \, dy \;\; \text{for} \;\; x \notin supp\, f,
\end{equation}
where kernel $K$ satisfies the following conditions:\\
standard size condition:
\begin{equation}
\label{formula 2 cond} \mdl{K(x,y)} \; \leq \; \frac{1}{\mdl{x-y}}
\end{equation}
and perfect cancelation condition:
\begin{equation}
\label{formula 3 cond}\mdl{K(x,y)-K(x,y^\prime)} +
\mdl{K(x,y)-K(x^\prime,y)} \; = \; 0
\end{equation}
whenever $x,x^\prime \in I \in D$, $y,y^\prime \in J \in D$, $I \cap
J = \varnothing$.

Note that, by perfect dyadic cancellation condition, $K(x,y)$ is
constant on $I^+ \times I^-$ and $I^- \times I^+$ for any dyadic interval $I \in D$.
We define
$$
K^+_I \; := \; K(x,y), \;\;\;\;\; x\in I^+, \;\; y\in I^-
$$
and
$$
K^-_I \; := \; K(x,y), \;\;\;\;\; x\in I^-, \;\; y\in I^+.
$$
Then
$$
K(x,y) = \sum_{I\in D} K^+_I \chi_{I^+\times I^-} + K^-_I
\chi_{I^-\times I^+}.
$$

Then perfect dyadic singular operator can be written as
\begin{equation}
\label{formula 4 Tf(x)}Tf(x) \; = \; \sum_{I\in D} K^+_I f_{I^-}
\chi_{I^+}(x) + K^-_I f_{I^+} \chi_{I^-}(x),
\end{equation}
where $f_I := \int_I f(x) dx$.

Now we are going to rewrite $T$ in a more convenient form. It is
easy to see that
\begin{eqnarray}
\label{formula 5 K+K}\nonumber && K^+_I f_{I^-} \chi_{I^+}(x) +
K^-_I f_{I^+}
\chi_{I^-}(x)\\
\nonumber &=& \frac{K^+_I f_{I^-} + K^-_I f_{I^+}}{2}\left(
\chi_{I^+}(x) + \chi_{I^-}(x) \right) + \frac{K^+_I f_{I^-} - K^-_I
f_{I^+}}{2}\left( \chi_{I^+}(x) - \chi_{I^-}(x) \right)\\
&=& \frac{1}{2}\left( K_I^+ f_{I^-} + K_I^- f_{I^+} \right)
\chi_I(x) + \frac{1}{2}\left( K_I^+ f_{I^-} - K_I^- f_{I^+} \right)
\sqrt{\mdl{I}} \; h_I(x),
\end{eqnarray}
where $\{h_I\}_{I\in D}$ is the Haar system of functions, normalized
in $L^2$, $h_I := \mdl{I}^{-1/2}(\chi_{I^+} - \chi_{I^-})$. Thus,
$Tf(x)$ can be written as:
\begin{equation}
\label{formula 6 Tf(x)}Tf(x) = \sum_{I\in D} \frac{1}{2} \left( K_I^+
f_{I^-} + K_I^- f_{I^+} \right) \chi_I + \frac{1}{2} \left( K_I^+ f_{I^-} - K_I^- f_{I^+} \right) h_I \sqrt{\mdl{I}}.
\end{equation}
Now let us handle the coefficients $\left( K_I^+ f_{I^-} + K_I^-
f_{I^+} \right)$ and $\left( K_I^+ f_{I^-} - K_I^- f_{I^+} \right)$:
\begin{eqnarray}
\label{formula 7 K+K}\nonumber && K_I^+ f_{I^-} + K_I^- f_{I^+} \; =
\; K_I^+\int_{I^-}f + K_I^-\int_{I^+}f\\
\nonumber &=& \frac{1}{2}\left( K_I^+ + K_I^- \right)\left(
\int_{I^-}f + \int_{I^+}f \right) - \frac{1}{2}\left( K_I^+ - K_I^-
\right)\left( \int_{I^+}f - \int_{I^-}f \right)\\
& = & \frac{1}{2}\left( K_I^+ + K_I^- \right) \mdl{I} \, m_If -
\frac{1}{2}\left( K_I^+ - K_I^- \right)\sqrt{\mdl{I}}\left\langle
f;h_I \right\rangle
\end{eqnarray}
and, similarly, replacing $K_I^-$ by $-K_I^-$,
\begin{equation}
\label{formula 8 K-K}K_I^+ f_{I^-} - K_I^- f_{I^+} \; = \;
\frac{1}{2}\left( K_I^+ - K_I^- \right) \mdl{I} \, m_If -
\frac{1}{2}\left( K_I^+ + K_I^- \right)\sqrt{\mdl{I}}\left\langle
f;h_I \right\rangle.
\end{equation}
We plug (\ref{formula 7 K+K}) and (\ref{formula 8 K-K}) in
(\ref{formula 6 Tf(x)}) and obtain a formal representation of the perfect dyadic operator $T$:
\begin{eqnarray}
\label{formula 9 Tf(x)}\nonumber Tf(x) &=& \frac{1}{4}\sum_{I\in
D}\left( K_I^+ + K_I^- \right) \mdl{I} \, m_If \, \chi_I -
\left( K_I^+ - K_I^- \right)\sqrt{\mdl{I}}
\left\langle f;h_I \right\rangle \chi_I\\
&& + \left( K_I^+ - K_I^- \right)
\mdl{I}^{3/2} \, m_If \, h_I - \left( K_I^+
+ K_I^- \right)\mdl{I}\left\langle f;h_I \right\rangle h_I.
\end{eqnarray}

We can represent $T^*$ in a similar way:
\begin{eqnarray}
\label{formula 10 T*f(x)}\nonumber T^*f(x) &=& \sum_{I\in D}
K_I^-f_{I^-}\chi_{I^+}(x) + K_I^+f_{I^+}\chi_{I^-}(x)\\
\nonumber &=& \sum_{I\in D} \frac{1}{2}\left( K_I^-f_{I^-} +
K_I^+f_{I^+} \right) \chi_I(x) + \frac{1}{2}\left( K_I^-f_{I^-} -
K_I^+f_{I^+} \right)\sqrt{\mdl{I}} \, h_I(x)\\
\nonumber &=& \frac{1}{4}\sum_{I\in D}\left( K_I^+ + K_I^- \right)
\mdl{I} \,m_If\, \chi_I(x) + \left( K_I^+ -
K_I^- \right)\sqrt{\mdl{I}} \left\langle f;h_I \right\rangle
\chi_I(x)\\
&& - \left( K_I^+ - K_I^- \right)
\mdl{I}^{3/2} m_If \, h_I(x) - \left( K_I^+
+ K_I^- \right) \mdl{I} \left\langle f;h_I \right\rangle h_I(x).
\end{eqnarray}


\section{Decomposition of $T$, additional assumptions.}

In this section we prove Theorem \ref{t:T(1)} modulo the unweighted boundedness of the operator $T_1$, which is done in a more general weighted case in the next section.

We assume, in addition to $T$ being perfect dyadic, that $T(1)$ and $T^*(1)$ are both in the dyadic $BMO^d$, with the $BMO^d$ norm of at most $Q$. We will also assume that $T$ satisfies testing conditions $\langle T h_I, h_I \rangle \leqslant Q$ for all dyadic intervals $I$.


\subsection{The $BMO^d$ condition.} A function $b(x)$, or a sequence $\{b_I\}_{I\in D} = \{\langle b; h_I \rangle\}_{I\in D}$ is in the dyadic $BMO^d$ whenever the sequence $b_I^2$ is a Carleson sequence, i.e. there is a finite constant C such that for all dyadic intervals $J\in D$ we have 
$$
\frac{1}{|J|} \sum_{I\in D(J)} b_{I}^2 \leqslant C.
$$
The best such constant C is called the dyadic $BMO^d$ norm of $b$.

Let $T(1)$, $T^*(1) \in BMO^d$. We want to show that
$$
\frac{1}{|J|} \sum_{I\in D(J)} (K_I^+ - K_I^-)^2 |I|^3 \leqslant 16Q.
$$

 First let us plug
$f(x)=1(x)$ in $T$ and $T^*$:
$$
T(1) \; = \; \frac{1}{4}\sum_{I\in D}\left( K_I^+ + K_I^- \right)
\mdl{I}\chi_I(x) + \left( K_I^+ - K_I^-
\right) \mdl{I}^{3/2}h_I(x)
$$
and
$$
T^*(1) \; = \; \frac{1}{4}\sum_{I\in D}\left( K_I^+ + K_I^- \right)
\mdl{I}\chi_I(x) - \left( K_I^+ - K_I^-
\right) \mdl{I}^{3/2}h_I(x).
$$

A function $f(x)$ belongs to the dyadic class $BMO^d$ whenever
$\sup_{J\in D}\frac{1}{\mdl{J}}\sum_{I\in D(J)}\langle f;h_I
\rangle^2 \leq Q$, the smallest such constant is the $BMO^d$-norm of
the function $f(x)$. So, let us find $\left\langle T(1);h_J
\right\rangle$ and $\left\langle T^*(1);h_J \right\rangle$:
\begin{eqnarray}
\label{formula 11 T(1);h}\nonumber \left\langle T(1);h_J
\right\rangle &=& \frac{1}{4}\sum_{I\in D}\left( K_I^+ + K_I^-
\right) \mdl{I}\left\langle \chi_I(x);h_J(x) \right\rangle +
\left( K_I^+ - K_I^- \right)
\mdl{I}^{3/2}\left\langle h_I(x);h_J(x) \right\rangle\\
\nonumber &=& \frac{1}{4}\sum_{I\in D(J^+)}\left( K_I^+ + K_I^-
\right) \mdl{I}\frac{\mdl{I}}{\sqrt{\mdl{J}}} -
\frac{1}{4}\sum_{I\in D(J^-)}\left( K_I^+ + K_I^- \right)\mdl{I}
\frac{\mdl{I}}{\sqrt{\mdl{J}}}\\
&& + \frac{1}{4}\left( K_J^+ - K_J^- \right) \mdl{J}^{3/2}.
\end{eqnarray}
Similarly,
\begin{eqnarray}
\label{formula 12 T*(1);h}\nonumber \left\langle T^*(1);h_J
\right\rangle &=& \frac{1}{4}\sum_{I\in D}\left( K_I^+ + K_I^-
\right) \mdl{I}\left\langle \chi_I(x);h_J(x) \right\rangle -
\left( K_I^+ - K_I^- \right)
\mdl{I}^{3/2}\left\langle h_I(x);h_J(x) \right\rangle\\
\nonumber &=& \frac{1}{4}\sum_{I\in D(J^+)}\left( K_I^+ + K_I^-
\right) \mdl{I}\frac{\mdl{I}}{\sqrt{\mdl{J}}} -
\frac{1}{4}\sum_{I\in D(J^-)}\left( K_I^+ + K_I^- \right)\mdl{I}
\frac{\mdl{I}}{\sqrt{\mdl{J}}}\\
&& - \frac{1}{4}\left( K_J^+ - K_J^- \right) \mdl{J}^{3/2}.
\end{eqnarray}
Let
$$
\alpha_J \; := \; \sum_{I\in D(J^+)}\left( K_I^+ + K_I^- \right)
\frac{\mdl{I}^2}{\sqrt{\mdl{J}}} - \sum_{I\in D(J^-)}\left( K_I^+ +
K_I^- \right) \frac{\mdl{I}^2}{\sqrt{\mdl{J}}}
$$
and
$$
\beta_J \; := \; \left( K_J^+ - K_J^- \right)\mdl{J}^{3/2},
$$
then
$$
\left\langle T(1);h_J \right\rangle \; = \; \frac{1}{4}\left( \alpha_J +
\beta_J \right), \;\;\;\;\; \left\langle T^*(1);h_J \right\rangle \; =
\; \frac{1}{4}\left( \alpha_J - \beta_J \right).
$$
If $T(1)\in BMO^d$ with $\norm{T(1)}{BMO^d} \leq Q$, then for every
dyadic interval $J\in D$
\begin{equation}
\label{formula 13 <Q} \frac{1}{\mdl{J}}\sum_{I\in D(J)}\left\langle
T(1);h_I \right\rangle^2 \; = \; \frac{1}{16\mdl{J}}\sum_{I\in D(J)}
\left( \alpha_I + \beta_I \right)^2 \; \leq \; Q.
\end{equation}
Similarly, $T^*(1)\in BMO^d$ with $\norm{T^*(1)}{BMO^d}\leq Q$
implies that
\begin{equation}
\label{formula 14 <Q} \forall J\in D \;\;\;\;\;
\frac{1}{\mdl{J}}\sum_{I\in D(J)} \left\langle T^*(1);h_I
\right\rangle^2 \; = \; \frac{1}{16\mdl{J}}\sum_{I\in D(J)} \left(
\alpha_I - \beta_I \right)^2 \; \leq \; Q.
\end{equation}
Let us see how conditions (\ref{formula 13 <Q}) and (\ref{formula 14
<Q}) imply that
$$
\forall J\in D \;\;\;\;\; \frac{1}{\mdl{J}}\sum_{I\in D(J)} \beta_I^2 \;
= \; \frac{1}{\mdl{J}}\sum_{I\in D}\left( K_I^+ - K_I^- \right)^2
\mdl{I}^3 \; \leq \; 16Q.
$$
The sum $\frac{1}{\mdl{J}}\sum_{I\in D(J)}\beta_I^2$ can be written as:
\begin{eqnarray*}
\frac{1}{\mdl{J}}\sum_{I\in D(J)}\beta_I^2 &=&
\frac{1}{\mdl{J}}\sum_{I\in D(J)}\left( \frac{\beta_I}{2} -
\frac{\alpha_I}{2} + \frac{\beta_I}{2} + \frac{\alpha_I}{2} \right)^2\\
&=& \frac{1}{\mdl{J}}\sum_{I\in D(J)}\left( \frac{\beta_I}{2} -
\frac{\alpha_I}{2} \right)^2 + \frac{2}{\mdl{J}}\sum_{I\in D(J)}\left(
\frac{\beta_I}{2} - \frac{\alpha_I}{2} \right)\left( \frac{\beta_I}{2} +
\frac{\alpha_I}{2} \right) + \frac{1}{\mdl{J}}\sum_{I\in D(J)}\left(
\frac{\beta_I}{2} + \frac{\alpha_I}{2} \right)^2\\
&\leq& 8Q + \frac{1}{2\mdl{J}}\sum_{I\in D(J)}\left( \beta_I - \alpha_I
\right) \left( \beta_I + \alpha_I \right).
\end{eqnarray*}
By Cauchy-Schwarz
$$
\frac{1}{\mdl{J}}\sum_{I\in D(J)}\left( \beta_I - \alpha_I \right) \left( \beta_I
+ \alpha_I \right) \; \leq \left( \frac{1}{\mdl{J}}\sum_{I\in D(J)}\left(
\beta_I - \alpha_I \right)^2 \right)^{\frac{1}{2}} \left(
\frac{1}{\mdl{J}}\sum_{I\in D(J)}\left( \beta_I + \alpha_I \right)^2
\right)^{\frac{1}{2}} \leq \; 16Q.
$$
Hence, for every dyadic interval $J\in D$
\begin{equation}
\label{formula 15 <9Q} \frac{1}{\mdl{J}}\sum_{I\in D(J)} \beta_I^2 \; =
\; \frac{1}{\mdl{J}}\sum_{I\in D}\left( K_I^+ - K_I^- \right)^2
\mdl{I}^3 \; \leq \; 16Q.
\end{equation}


\subsection{Testing conditions.}
We assume that $T$ satisfies testing conditions
\begin{equation}
\label{formula 16 Th;h<Q} \forall J\in D \;\;\;\;\; \left\langle
Th_J;h_J \right\rangle \; \leq \; Q.
\end{equation}
Let us write $\left\langle Th_J;h_J \right\rangle$ as follows:
\begin{eqnarray*}
\mdl{\left\langle Th_J;h_J \right\rangle} &=& \frac{1}{4}\left|
\sum_{I\in D}\left( K_I^+ + K_I^- \right) \left\langle h_J;\chi_I
\right\rangle \left\langle \chi_I;h_J \right\rangle - \left( K_I^+ - K_I^- \right)\sqrt{\mdl{I}} \left\langle h_J;h_I
\right\rangle \left\langle \chi_I;h_J \right\rangle \right.\\
&& \left. + \left( K_I^+ - K_I^- \right)\sqrt{\mdl{I}}
\left\langle h_J;\chi_I \right\rangle \left\langle h_I;h_J
\right\rangle - \left( K_I^+ + K_I^- \right)\mdl{I}
\left\langle h_J;h_I \right\rangle \left\langle h_I;h_J
\right\rangle \right|\\
&=& \left| \frac{1}{4}\sum_{I\in D(J^+)}\left( K_I^+ + K_I^- \right)
\left( \frac{\mdl{I}}{\sqrt{\mdl{J}}} \right)^2 +
\frac{1}{4}\sum_{I\in D(J^-)}\left( K_I^+ + K_I^- \right) \left(
-\frac{\mdl{I}}{\sqrt{\mdl{J}}} \right)^2\right.\\
&& \left. - \frac{1}{4}\left( K_J^+ + K_J^- \right)\mdl{J} \right|\\
&=& \mdl{\frac{1}{4\mdl{J}}\sum_{I\in D(J)}\left( K_I^+ + K_I^-
\right)\mdl{I}^2 - \frac{1}{2}\left( K_J^+ + K_J^- \right) \mdl{J}}
\; \leq \; Q.
\end{eqnarray*}
Note that by standard size condition on the kernel $K$, we know that
$\mdl{K_J^+ + K_J^-}\mdl{J} \leq 2Q$, so
\begin{equation}
\label{formula 17 <8Q} \mdl{\frac{1}{\mdl{J}}\sum_{I\in D(J)}\left(
K_I^+ + K_I^- \right)\mdl{I}^2} \; \leq \; 8Q.
\end{equation}
So, a perfect dyadic Calder\'{o}n-Zygmund singular integral operator
$T$ satisfying the $T(1)$ conditions $\|T(1)\|_{BMO^2}, \|T^*(1)\|_{BMO^d} \leqslant Q$ and the dyadic testing conditions $\left( \left\langle
Th_J;h_J \right\rangle \leq Q \;\;\; \forall J\in D \right)$, can be
written as:
\begin{eqnarray}
\label{formula 18 Tf(x)} \nonumber Tf(x) &=& \sum_{I\in D}
K_I^+f_{I^-}\chi_{I^+}(x) + K_I^-f_{I^+}\chi_{I^-}(x)\\
\nonumber &=& \frac{1}{4}\sum_{I\in D}\left( K_I^+ + K_I^- \right)
\mdl{I} m_If\, \chi_I(x) - \left( K_I^+ -
K_I^- \right)\sqrt{\mdl{I}}\left\langle f;h_I \right\rangle
\chi_I(x)\\
&& +  \left( K_I^+ - K_I^- \right)
\mdl{I}^{3/2} m_If \, h_I(x) - \left( K_I^+
+ K_I^- \right) \mdl{I} \left\langle f;h_I \right\rangle h_I(x).
\end{eqnarray}
Moreover, coefficients $K_I^+$ and $K_I^-$ satisfy the following
conditions:\\
size conditions
\begin{equation}
\label{formula 19 size cond} \forall J\in D \;\;\;\;\; \mdl{K_J^+ +
K_J^-} \mdl{I} \; \leq \; 2 Q,
\end{equation}
$T(1)$ conditions
\begin{equation}
\label{formula 20 T(1) cond} \forall J\in D \;\;\;\;\;
\frac{1}{\mdl{J}}\sum_{I\in D(J)}\left( K_I^+ - K_I^- \right)^2
\mdl{I}^3 \; \leq \; 16 Q
\end{equation}
and testing conditions
\begin{equation}
\label{formula 21 test cond} \forall J\in D \;\;\;\;\;
\mdl{\frac{1}{\mdl{J}}\sum_{I\in D(J)}\left( K_I^+ + K_I^-
\right)\mdl{I}^2} \; \leq \; 8 Q.
\end{equation}
By (\ref{formula 18 Tf(x)}) operator $T$ can be written as 

\begin{equation} \label{decomposition} 
T = \frac{1}{4}
\left( T_1 - T_2 + T_3 - T_4 \right),
\end{equation}
where
\begin{equation}
\label{formula 22 T1} T_1f(x) \; = \; \sum_{I\in D}\left( K_I^+ +
K_I^- \right) \mdl{I} m_If\, \chi_I(x),
\end{equation}
\begin{equation}
\label{formula 23 T2} T_2f(x) \; = \; \sum_{I\in D}\left( K_I^+ -
K_I^- \right)\sqrt{\mdl{I}}\left\langle f;h_I \right\rangle
\chi_I(x),
\end{equation}
\begin{equation}
\label{formula 24 T3} T_3f(x) \; = \; \sum_{I\in D} \left( K_I^+ -
K_I^- \right) \mdl{I}^{3/2} m_If \, h_I(x),
\end{equation}
\begin{equation}
\label{formula 25 T4} T_4f(x) \; = \; \sum_{I\in D}\left( K_I^+ +
K_I^- \right) \mdl{I} \left\langle f;h_I \right\rangle h_I(x),
\end{equation}
since all four operators are known well-defined linear operators bounded on $L^2$. Operators $T_3$ and $T_2$ are dyadic paraproduct and its adjoint, that are both well defined and bounded in $L^2$ since the sequence $b_I = (K_I^+ - K_I^-) |I|^{3/2}$ satisfies $\sum_{I\in D(J)} b_I^2 \leqslant C |J|$ for any $J \in D$ by (\ref{formula 20 T(1) cond}). $T_4$ is a martingale transform with symbol bounded uniformly by (\ref{formula 19 size cond}). In the next section we show that $T_1$ is bounded on $L^2$ and on all weighted spaces for weights in $A_2^d$.

This completes the proof of Theorem \ref{t:T(1)}.


\section{Weighted $T(1)$ theorem for perfect dyadic operators}

Let $w$ be a weight, i.e. almost everywhere positive locally integrable function on the real line. Let a weight $w$ be such that $w^{-1}$ is a weight as well. Assume that $w$ is in the dyadic Muckenhoupt class $A_2^d$, i.e. 
$$
[w]_{A_2^d} := \sup_{I\in D} m_I w \, m_I (w^{-1}) \; <\infty.
$$ 

In order to show that the weighted $L^2(w)$ norm of the operator $T$
depends on the $A_2$-constant of the weight $w$ at most linearly, it
is enough to show the linear bound for each of the operators $T_i$,
$i=1,2,3,4$.

First let us show that $T_1$ is bounded on $L^2(w)$ and its norm
depends on the $A_2$-constant of the weight $w$ at most linearly. We
will do it by duality, we will show that $\forall f\in L^2(w)$ and
$\forall g\in L^2(w^{-1})$
$$
\left\langle T_1 f;g \right\rangle \; = \; \sum_{I\in D} \left( K_I^+ +
K_I^- \right) \mdl{I}^2 m_If \, m_Ig \; \leq \; C
[{w}]_{A_2^d}\norm{f}{L^2(w)}\norm{g}{L^2(w^{-1})}
$$
or, alternatively, $\forall f,g\in L^2$
\begin{equation}
\label{formula 26 linear T1} \sum_{I\in D} \left( K_I^+ + K_I^-
\right) \mdl{I}^2 m_I\left( fw^{-1/2} \right)\, m_I\left( gw^{1/2}
\right) \; \leq \; C [{w}]_{A_2^d} \norm{f}{L^2} \norm{g}{L^2}.
\end{equation}
Without loss of generality we may assume that coefficients $(K_I^+ + K_I^-)$ are all non-negative. We are going to use the following version of the bilinear embedding
theorem from \cite{NazarovTreilVolberg:99}.
\begin{theorem}[Nazarov, Treil, Volberg]
\label{theorem 1 NTrVo2} Let $v$ and $w$ be weights. Let $\left\{
a_I \right\}$ be a sequence of nonnegative numbers, s.t. for
all dyadic intervals $J\in D$ the following three inequalities hold
with some constant $Q>0$:
\begin{equation}
\label{formula 27 theorem1} \frac{1}{\mdl{J}}\sum_{I\in
D(J)}a_I\, m_Iw \, m_Iv \; \leq \; Q,
\end{equation}
\begin{equation}
\label{formula 28 theorem1} \frac{1}{\mdl{J}}\sum_{I\in
D(J)}a_I\, m_Iw \; \leq \; Q \, m_Iw,
\end{equation}
\begin{equation}
\label{formula 29 theorem1} \frac{1}{\mdl{J}}\sum_{I\in
D(J)}a_I\, m_Iv \; \leq \; Q\, m_Iv,
\end{equation}
then for any two nonnegative functions $f,g\in L^2$
\begin{equation}
\label{formula 30 theorem1}\sum_{I\in D} a_I\, m_I\left(
fv^{1/2} \right) \, m_I\left( gw^{1/2} \right) \; \leq \; C\, Q\,
\norm{f}{L^2} \norm{g}{L^2}.
\end{equation}
\end{theorem}


Applying this theorem with $v=w^{-1}$ and $a_I = 
\left( K_I^+ + K_I^- \right)\mdl{I}^2$, we can see that the desired
bound (\ref{formula 26 linear T1}) will hold if we can prove the
following three inequalities for every dyadic interval $J\in D$:
\begin{equation}
\label{formula 31 ineq1of3}\frac{1}{\mdl{J}}\sum_{I\in D(J)}\left(
K_I^+ + K_I^- \right)\mdl{I}^2 m_Iw \; m_I(w^{-1}) \; \leq \;
C[{w}]_{A_2^d},
\end{equation}
\begin{equation}
\label{formula 32 ineq2of3}\frac{1}{\mdl{J}}\sum_{I\in D(J)}\left(
K_I^+ + K_I^- \right)\mdl{I}^2 m_Iw \; \leq \; C[w]_{A_2^d} m_Jw,
\end{equation}
\begin{equation}
\label{formula 33 ineq3of3}\frac{1}{\mdl{J}}\sum_{I\in D(J)}\left(
K_I^+ + K_I^- \right)\mdl{I}^2 m_I(w^{-1}) \; \leq \; C[w]_{A_2^d}
m_J(w^{-1}).
\end{equation}
It turns out that all these inequalities follow from the testing
conditions (\ref{formula 21 test cond}). It is easy to see that
inequality (\ref{formula 31 ineq1of3}) follows from (\ref{formula 21
test cond}) and the definition of $A_2$-constant $\left( m_Iw \,
m_I(w^{-1}) \leq [w]_{A_2^d} \right)$.

In order to see that inequalities (\ref{formula 32 ineq2of3}) and
(\ref{formula 33 ineq3of3}) are true, we will need the following
lemma, which can be found, for example, in \cite{Beznosova:08}:
\begin{lemma}
\label{lemma1 Be} Let $v$ be a weight, such that $v^{-1}$ is a
weight as well, and let $\left\{ \lambda_I \right\}$ be a Carleson
sequence of nonnegative numbers, that is, there exists a constant
$Q>0$ s.t.
$$
\forall J\in D \;\;\;\;\; \frac{1}{\mdl{J}}\sum_{I\in D(J)}
\lambda_I \; \leq \; Q,
$$
then
$$
\forall J\in D \;\;\;\;\; \frac{1}{\mdl{J}}\sum_{I\in
D(J)}\frac{\lambda_I}{m_I(v^{-1})} \; \leq \; 4Q\, m_Jv
$$
and, therefore, if $v\in A_2$ then for any $J\in D$ we have
$$
\frac{1}{\mdl{J}}\sum_{I\in D(J)} m_Iv \, \lambda_I \; \leq \; 4Q\,
\norm{v}{A_2^d} \, m_Jv.
$$
\end{lemma}
Applying this lemma to $\lambda_I = \left( K_I^+ + K_I^- \right)
\mdl{I}^2$ and $v=w$, we obtain that (\ref{formula 32 ineq2of3})
follows from (\ref{formula 21 test cond}). If we take $v=w^{-1}$ and
observe that $\norm{w^{-1}}{A_2^d} = \norm{w}{A_2^d}$, we can see
that (\ref{formula 33 ineq3of3}) follows from (\ref{formula 21 test
cond}) as well.

So,
$$
\norm{T_1}{L^2(w)\rightarrow L^2(w)} \; \leq \; C \norm{w}{A_2^d}.
$$
Let us analyze operators $T_3$ and $T_2$ now.
$$
T_3f(x) \; = \; \sum_{I\in D} \left( K_I^+ - K_I^- \right)
\mdl{I}^{3/2} m_If \, h_I(x),
$$
$$
T_2f(x) \; = \; \sum_{I\in D}\left( K_I^+ - K_I^-
\right)\sqrt{\mdl{I}}\left\langle f;h_I \right\rangle \chi_I(x),
$$
We first note that $T_3$ is a paraproduct operator and $T_2$ is its
adjoint:
$$
T_3f(x) \; = \; \pi_b f(x) \; = \; \sum_{I\in D} m_If \; b_I \;
h_I(x),
$$
$$
T_2f(x) \; = \; \pi_b^* f(x) \; = \; \sum_{I\in D} b_I \left\langle
f;h_I \right\rangle \frac{1}{\mdl{I}} \chi_I(x)
$$
with the sequence $b_I = \left( K_I^+ - K_I^- \right)\mdl{I}^{3/2}$.

In order for dyadic paraproduct and its adjoint to be bounded in
$L^2$ we need $b$ to be in the $BMO^d$, i.e. we need $\left\{ b_I^2
\right\}_{I\in D}$ to be a Carleson sequence, which is
$$
\forall J\in D \;\;\;\;\; \frac{1}{\mdl{J}}\sum_{I\in D}\left( K_I^+
- K_I^- \right)^2\mdl{I}^{3} \; \leq \; Q
$$
and coincides with $T(1)$ conditions (\ref{formula 20 T(1) cond}).

It has been shown in \cite{Beznosova:08} that the norm of the dyadic
paraproduct (and, by symmetry of $A_2$-constant, of its adjoint as
well) depends on the $A_2$-constant of the weight $w$ at most
linearly, i.e.
$$
\norm{T_3}{L^2(w)\rightarrow L^2(w)} + \norm{T_2}{L^2(w)\rightarrow
L^2(w)} \; \leq \; C\, Q^{1/2} [w]_{A_2^d}.
$$

And, finally, the linear bound on the operator $T_4$,
$$
T_4f(x) \; = \; \sum_{I\in D}\left( K_I^+ + K_I^- \right) \mdl{I}
\left\langle f;h_I \right\rangle h_I(x),
$$
is implied by the uniform boundedness of the symbol (i.e. size conditions (\ref{formula 19 size cond})) and
Wittwer's weighted linear bound on the martingale transform (see \cite{Wittwer:00}).

Since all parts of $T$ obey linear bounds on their $L^2(w)$ norms
with respect to the $A_2$-constant of the weight $w$, bound on the
norm of $T$ is linear as well,
$$
\norm{T}{L^2(w)\rightarrow L^2(w)} \; \leq \; C \norm{w}{A_2},
$$
where constant $C$ depends on the constant in the size condition
(\ref{formula 2 cond}), dyadic $BMO^d$-norms of $T(1)$ and $T^*(1)$
(at most linearly) and the constant in testing conditions
(\ref{formula 16 Th;h<Q}) (at most linearly as well). Therefore, we have the following theorem.

\begin{theorem} \label{t:A2conjecture}
Let $T$ be a perfect dyadic operator on $\R$ such that 
$$\|T(1)\|_{BMO^d} \leqslant Q \;\;\; \text{and} \;\;\; \|T^*(1)\|_{BMO^d} \leqslant Q $$
and
$$ \forall I\in D \;\;\; \langle T h_I , h_I \rangle \leqslant Q $$.
Then $T$ is bounded on $L^2(w)$ and 
$$
\|T f\|_{L^2(w)} \leqslant C Q \|f\|_{L^2(w)},
$$
with some constant $C$ independent of the operator $T$.
\end{theorem}

Theorem \ref{t:A2conjecture} in the case $w(x)=1$ is a $T(1)$ theorem, so we can view it as a weighted version of the $T(1)$ theorem. It is also interesting that the $L^2(w)$ norm of the operator $T$ after proper normalization (we assume that the decay constant of the kernel is $1$) only depends the $BMO^d$ norms of $T(1)$ and $T^*(1)$ and the constant in testing conditions. It also depends on these constants at most linearly.


\section{Two weight boundedness of perfect dyadic operators.}

We start with the decomposition (\ref{decomposition}), 
$$
T f(x) = \frac{1}{4} (T_1 f(x) + T_2 f(x) + T_3 f(x) + T_4 f(x)).
$$
First we consider 
$$
T_1 f = \sum_{I \in D} (K_I^+ + K_I^-) |I| m_If \chi_I(x).
$$
Similarly to the one weight case, by duality using Theorem \ref{theorem 1 NTrVo2}, we know it is bounded from $L^2(v)$ in $L^2(u)$ (with the norm bounded by $C Q_1$) whenever for every dyadic interval $J \in D$ the following three conditions hold simultaneously
\begin{equation}\label{e:T1-1}
\frac{1}{|J|} \sum_{I\in D(J)} (K_I^+ + K_I^-) |I|^2 m_I(v^{-1}) m_I u \leqslant Q_1,
\end{equation}
\begin{equation}\label{e:T1-2}
\frac{1}{|J|} \sum_{I\in D(J)} (K_I^+ + K_I^-) |I|^2  m_I u \leqslant Q_1 m_J u,
\end{equation}
\begin{equation}\label{e:T1-3}
\frac{1}{|J|} \sum_{I\in D(J)} (K_I^+ + K_I^-) |I|^2 m_I(v^{-1}) \leqslant Q_1 m_J (v^{-1}).
\end{equation}

 Second, consider the paraproduct $T_3$ and its adjoint $T_2$. In \cite{BCMP} we have sufficient conditions for the two weight boundedness for the dyadic paraproduct:

\begin{theorem}\label{T:para1} Let $\pi_b$ be the dyadic paraproduct operator associated to the sequence $b = \{b_I\}_{I\in D}$:
$$
\pi_bf := \sum_{I\in {D}} m_If \,b_I\, h_I.
$$
Let $(v,u)$ be a pair of measurable functions on $\R$ such that $u$ and
$v^{-1}$, the reciprocal of $v$, are 
weights  on $\R$ and such that

(i) $(v,u)\in {A}_2^d$, that is
$ \; [v,u]_{{A}_2^d} := \sup_{I\in {D}}m_I (v^{-1})\, m_I u <\infty.$

(ii) 
Assume that there is a constant $\mathcal{C}_{v,u}>0$ such that
\[ \sum_{I\in {D}(J)} |\Delta_I u|^2|I|\, m_I(v^{-1}) \leq \mathcal{C}_{v,u} u(J) \quad\quad \mbox{for all} \; J\in {D},\]
where $\Delta_I u :=  m_{I_+}u-m_{I_-}u$, and $I_{\pm}$ are the right and left children of $I$.

Assume that $b\in Carl_{v,u}$, i.e.
there is a constant $\mathcal{B}_{v,u}>0$ such that
\[ \sum_{I\in {D}(J)} \frac{|b_I|^2}{m_Iu} \leq \mathcal{B}_{v,u} v^{-1}(J)\quad\quad \mbox{for all} \; J\in {D}.\]

Then $\pi_b$, the dyadic paraproduct associated to $b$, is bounded from
$L^2(v)$ in $L^2(u)$. 
Moreover,  there exist $C>0$ such that
$$\|\pi_bf\|_{L^2(u)}\leq C\sqrt{[v,u]_{{A}_2^d}\mathcal{B}_{v,u}}
\Big(\sqrt{[v,u]_{{A}_2^d}}+\sqrt{\mathcal{C}_{v,u}}\,\Big)\|f\|_{L^2(v)}\,.$$
\end{theorem}

Therefore, in order to be able to bound the paraproduct $T_3$ that has symbol $b_I = (K_I^+ - K_I^-) |I|^{3/2}$, we need the following three conditions to hold:

\begin{equation}\label{e:T3-1}
 [v,u]_{{A}_2^d} := \sup_{I\in {D}}m_I (v^{-1})\, m_I u <\infty,
\end{equation}
\begin{equation}\label{e:T3-2}
\sum_{I\in {D}(J)} |\Delta_I u|^2|I|\, m_I(v^{-1}) \leq \mathcal{C}_{v,u} u(J) \quad\quad \mbox{for all} \; J\in {D},
\end{equation}
\begin{equation}\label{e:T3-3}
\sum_{I\in {D}(J)} \frac{(K_I^+ - K_I^-)^2 |I|^{3}}{m_Iu} \leq \mathcal{B}_{v,u} v^{-1}(J)\quad\quad \mbox{for all} \; J\in {D}.
\end{equation}

The operator $T_2$ is adjoint to $T_3$, $T_2 = T_3^*$, so $T_2$ is bounded from $L^2(v)$ in $L^2 (u)$ whenever the paraproduct $T_3$ is bounded from $L^2(u^{-1})$ in $L^2(v^{-1})$. Therefore, in order for $T_2$ to be bounded from $L^2(v)$ in $L^2 (u)$ we need the following three conditions to hold:

\begin{equation}\label{e:T2-1}
 [u^{-1},v^{-1}]_{{A}_2^d} := \sup_{I\in {D}}m_I (u)\, m_I (v^{-1}) <\infty,
\end{equation}
\begin{equation}\label{e:T2-2}
\sum_{I\in {D}(J)} |\Delta_I (v^{-1})|^2|I|\, m_I u \leq \mathcal{C}_{u^{-1},v^{-1}} v^{-1}(J) \quad\quad \mbox{for all} \; J\in {D},
\end{equation}
\begin{equation}\label{e:T2-3}
\sum_{I\in {D}(J)} \frac{(K_I^+ - K_I^-)^2 |I|^{3}}{m_I (v^{-1})} \leq \mathcal{B}_{u^{-1},v^{-1}} u(J)\quad\quad \mbox{for all} \; J\in {D}.
\end{equation}


Finally, consider $T_4$, the martingale transform with symbol $\sigma = (K_I^+ +
K_I^-) |I|$, which is uniformly bounded by (\ref{formula 19 size cond}),
$$
T_4f(x) \; = \; \sum_{I\in D}\left( K_I^+ +
K_I^- \right) \mdl{I} \left\langle f;h_I \right\rangle h_I(x).
$$
Necessary and sufficient conditions were obtained in \cite{NazarovTreilVolberg:99}:
\begin{theorem}\label{t:NTV1} Let $\sigma = \{\sigma_I\}_{I\in D}$ be a sequence of signs $\pm$.
The family of operators $T_\sigma$ of martingale transforms with symbol $\sigma$ is  uniformly bounded from $L^2(v)$ to $L^2(u)$ if and only if the following assertions hold simultaneously:
\begin{equation}\label{e:T4-1}
\forall J \in D \;\;\; m_J u \; m_J v^{-1} \leqslant C <\infty,
\end{equation}
\begin{equation}\label{e:T4-2}
\forall J \in D \;\;\; \frac{1}{|J|} \sum_{I \in D(J)} (\Delta_I u)^2 m_I (v^{-1}) |I| \leqslant C m_I u,
\end{equation}
\begin{equation}\label{e:T4-3}
\forall J \in D \;\;\; \frac{1}{|J|} \sum_{I \in D(J)} (\Delta_I (v^{-1}))^2 m_I u |I| \leqslant C m_I (v^{-1}),
\end{equation}

and an operator $T_0$, defined by $T_0 f := \sum_{I\in D} \left |\frac{\Delta_I (v^{-1}) \Delta_I u}{m_I(v^{-1}) m_I u}\right | m_I f \chi_I (x)$, is bounded from $L^2(v)$ to $L^2(u)$, or, equivalently

\begin{equation}\label{e:T4-4}
\forall J \in D \;\;\; \frac{1}{|J|} \int_J \left(\sum_{I \in D(J)}  \left|\frac{\Delta_I u \Delta_I (v^{-1})}{m_Iu } \right| \chi_I(x) \right)^2 u dx \leqslant C m_J (v^{-1})
\end{equation}
and
\begin{equation}\label{e:T4-5}
\forall J \in D \;\;\; \frac{1}{|J|} \int_J \left(\sum_{I \in D(J)}   \left|\frac{\Delta_I u \Delta_I (v^{-1})}{ m_I(v^{-1})}\right| \chi_I(x) \right)^2 v^{-1} dx \leqslant C m_J u.
\end{equation}
\end{theorem}

Now we need to put all conditions we obtained for the operators $T_{1,2,3,4}$, conditions (\ref{e:T1-1}), (\ref{e:T1-2}), (\ref{e:T1-3}), (\ref{e:T3-1}), (\ref{e:T3-2}), (\ref{e:T3-3}), (\ref{e:T2-1}), (\ref{e:T2-2}), (\ref{e:T2-3}), (\ref{e:T4-1})(\ref{e:T4-2})(\ref{e:T4-3})(\ref{e:T4-4})(\ref{e:T4-5}), together. 
Note that (\ref{e:T3-1}), (\ref{e:T2-1}) and (\ref{e:T4-1}) are all the same joint dyadic $A_2^d$ condition, while condition (\ref{e:T1-1}) follows from the joint $A_2^d$ and (\ref{formula 21 test cond}) under no additional assumptions.
Note  also that conditions (\ref{e:T3-2}) and (\ref{e:T4-2}) are the same, conditions (\ref{e:T2-2}) and (\ref{e:T4-3}) also coincide. Therefore, we obtain the following theorem.

\begin{theorem}\label{t:Tuv}
Let $(v,u)$ be a pair of measurable functions, such that $u$ and $v^{-1}$, the reciprocal of $v$, are weights on $\R$ such that
\begin{equation}\label{e:uv-1}
(v,u)\in A_2^d \,,\text{i.e.} \;\; [v,u]_{A_2^d} := \sup_{I\in D} m_I (v^{-1}) m_I u <\infty,
\end{equation}
\begin{equation}\label{e:uv-2}
\forall J \in D \;\;\; \frac{1}{|J|} \sum_{I \in D(J)} (\Delta_I u)^2 m_I (v^{-1}) |I| \leqslant C m_I u,
\end{equation}
\begin{equation}\label{e:uv-3}
\forall J \in D \;\;\; \frac{1}{|J|} \sum_{I \in D(J)} (\Delta_I (v^{-1}))^2 m_I u |I| \leqslant C m_I (v^{-1}),
\end{equation}
and operator $T_0$ is bounded from $L^2(v)$ in $L^2(u)$, i.e.
\begin{equation}\label{e:uv-4}
\forall J \in D \;\;\; \frac{1}{|J|} \int_J \left(\sum_{I \in D(J)}  \left|\frac{\Delta_I u \Delta_I (v^{-1})}{m_Iu } \right| \chi_I(x) \right)^2 u(x) dx \leqslant C m_J (v^{-1})
\end{equation}
and
\begin{equation}\label{e:uv-5}
\forall J \in D \;\;\; \frac{1}{|J|} \int_J \left(\sum_{I \in D(J)}   \left|\frac{\Delta_I u \Delta_I (v^{-1})}{ m_I(v^{-1})}\right| \chi_I(x) \right)^2 v^{-1}(x) dx \leqslant C m_J u.
\end{equation}
Let $T$ be a perfect dyadic operator with $\|T(1)\|_{BMO^d}, \|T^*(1)\|_{BMO^d} \leqslant Q$ and for every $I \in D$ we have $\langle T h_I; h_I \rangle \leqslant Q$.

Then $T$ is bounded from $L^2 (v)$ to $L^2(u)$ whenever for any dyadic interval $J \in D$

\begin{equation}\label{e:Tuv-1}
\frac{1}{|J|} \sum_{I\in D(J)} (K_I^+ + K_I^-) |I|^2  m_I u \leqslant C m_J u,
\end{equation}
\begin{equation}\label{e:Tuv-2}
\frac{1}{|J|} \sum_{I\in D(J)} (K_I^+ + K_I^-) |I|^2 m_I(v^{-1}) \leqslant C m_J (v^{-1}),
\end{equation}
\begin{equation}\label{e:Tuv-3}
\frac{1}{|J|}\sum_{I\in {D}(J)} \frac{(K_I^+ - K_I^-)^2 |I|^{3}}{m_Iu} \leq C m_J (v^{-1})
\end{equation}
and
\begin{equation}\label{e:Tuv-4}
\frac{1}{|J|}\sum_{I\in {D}(J)} \frac{(K_I^+ - K_I^-)^2 |I|^{3}}{m_I (v^{-1})} \leq C m_J u. 
\end{equation}
\end{theorem}


\section{Sufficient conditions for the two weight boundedness of perfect dyadic operators under additional assumptions that weighted are in $A_\infty^d$}

In this section we will assume that $v^{-1}$ and $u$ are $A^d_\infty$ weights and show that conditions of Theorem \ref{t:Tuv} can be reduced in this case. A weight $w$ belongs to the class $A_\infty^d$ whenever
$$
[w]_{A_\infty^d} := \sup_{I\in D} m_I w e^{-m_I(\ln w)} \,<\, \infty.
$$
From \cite{mydisser} we have the following lemma

\begin{lemma} \label{t:littleoo}
Let $w$ be a weight, such that $\log w \in L^1_{loc}$, let $\{\lambda_I\}_{I\in D}$ be a Carleson sequence:
$$
\exists Q >0 \;\;s.t.\;\; \forall J\in D \;\;\; \frac{1}{|J|} \sum_{I \in D(J)} \lambda_I \leqslant Q.
$$
Then for every dyadic interval $J \in D$
$$
\frac{1}{|J|} \sum_{I\in D(J)} e^{m_I (\log w)} \lambda_I \leqslant 4 Q m_J w
$$
and therefore, if $w \in A_\infty^d$ then for every dyadic interval $J\in D$ we have:
$$
\frac{1}{|J|} \sum_{I\in D(J)} m_Iw \;\lambda_I \leqslant 4 Q [w]_{A_\infty^d} m_I w.
$$
\end{lemma}

In particular, $u,v^{-1} \in A_\infty^d$ implies that $u$ and $v^{-1}$ are in $RH_1^d$ (see \cite{BeRez}), where $RH_1^d$ is defined as follows
$$
w \in RH_1^d \;\;\;\Longleftrightarrow \;\;\; [w]_{RH_1^d} := \sup_{I\in D} m_I \left(\frac{w}{m_I w}\log \frac{w}{m_I w}\right) \;<\,\infty.
$$
We also know from \cite{BeRez} that $[w]_{RH_1^d} \leqslant \log 16 [w]_{A_\infty^d}$ and we have the following theorem that first appeared in \cite{Bu} without the sharp constant and in \cite{BeRez} it was stated in stronger form and the constant was traced. Here we state the strong form from \cite{BeRez}.
\begin{theorem} \label{t:Buckley}
Let $w$ be almost everywhere positive locally integrable function on the real line and $J$ be any interval. Then up to a numerical constant
$$
m_J \left(w \log w\right) - m_J w \log m_J w \approx \frac{1}{|J|} \sum_{I\in D(J)} \left(\frac{\Delta_I w}{m_I w}\right)^2 m_I w |I|.
$$
In particular case when $w$ is a weight in $RH_1^d$, we have that for every dyadic interval $J \in D$
\begin{equation}
\frac{1}{|J|} \sum_{I\in D(J)} \left(\frac{\Delta_I w}{m_I w}\right)^2 m_I w |I| \leqslant C [w]_{RH_1^d} m_J w.
\end{equation}
\end{theorem}

Firstly, note that conditions (\ref{e:uv-2}) and (\ref{e:uv-3}) follow from the fact that $u$ and $v^{-1}$ are in $A_\infty^d$, Theorem \ref{t:Buckley} and the joint $A_2^d$ condition (\ref{e:uv-1}).

Secondly, conditions (\ref{e:Tuv-1}) and (\ref{e:Tuv-2}) by Lemma \ref{t:littleoo} follow from the fact that $u, v^{-1} \in A_\infty^d$ and \ref{formula 21 test cond}. Therefore, under additional assumptions that $u,v^{-1} \in A_\infty^d$, we can simplify Theorem \ref{t:Tuv} as follows.

\begin{theorem}\label{t:Tuv-oo}
Let $(v,u)$ be a pair of measurable functions, such that $u$ and $v^{-1}$, the reciprocal of $v$, are $A_\infty^d$ weights on $\R$ such that
\begin{equation}\label{e:uv-1}
(v,u)\in A_2^d \,,\text{i.e.} \;\; [v,u]_{A_2^d} := \sup_{I\in D} m_I (v^{-1}) m_I u <\infty,
\end{equation}
and operator $T_0$ is bounded from $L^2(v)$ in $L^2(u)$, i.e.
\begin{equation}\label{e:uv-4}
\forall J \in D \;\;\; \frac{1}{|J|} \int_J \left(\sum_{I \in D(J)}  \left|\frac{\Delta_I u \Delta_I (v^{-1})}{m_Iu } \right| \chi_I(x) \right)^2 u(x) dx \leqslant C m_J (v^{-1})
\end{equation}
and
\begin{equation}\label{e:uv-5}
\forall J \in D \;\;\; \frac{1}{|J|} \int_J \left(\sum_{I \in D(J)}   \left|\frac{\Delta_I u \Delta_I (v^{-1})}{ m_I(v^{-1})}\right| \chi_I(x) \right)^2 v^{-1}(x) dx \leqslant C m_J u.
\end{equation}
Let $T$ be a perfect dyadic operator with $\|T(1)\|_{BMO^d}, \|T^*(1)\|_{BMO^d} \leqslant Q$ and for every $I \in D$ we have $\langle T h_I; h_I \rangle \leqslant Q$.

Then $T$ is bounded from $L^2 (v)$ to $L^2(u)$ whenever for any dyadic interval $J \in D$

\begin{equation}\label{e:Tuv-3}
\frac{1}{|J|}\sum_{I\in {D}(J)} \frac{(K_I^+ - K_I^-)^2 |I|^{3}}{m_Iu} \leq C m_J (v^{-1})
\end{equation}
and
\begin{equation}\label{e:Tuv-4}
\frac{1}{|J|}\sum_{I\in {D}(J)} \frac{(K_I^+ - K_I^-)^2 |I|^{3}}{m_I (v^{-1})} \leq C m_J u. 
\end{equation}
\end{theorem}


\bibliographystyle{plain}

\end{document}